\RequirePackage{fix-cm}
\documentclass[smallextended]{svjour3}       
\smartqed  
\usepackage{geometry}
\usepackage{graphicx}
\usepackage{graphicx,amssymb,amsmath}
\usepackage{epsfig}
\usepackage{subfigure}
\usepackage{cite}
\usepackage{booktabs}
\usepackage{array}
\usepackage{color}
\usepackage{enumerate}
\usepackage{algorithm}
\usepackage{algorithmic}
\usepackage{float}
\usepackage{multirow}
\usepackage{rotating}
\usepackage[misc]{ifsym}
\usepackage{bm}
\usepackage{mathrsfs}
\definecolor{pred}{RGB}{148,55,61}
\definecolor{blue}{rgb}{0,0,0.9}
\definecolor{red}{rgb}{0.9,0,0}
\definecolor{green}{rgb}{0,0.9,0}

\newtheorem{assumption}{Assumption}
\numberwithin{equation}{section}

\begin{document}

\title{
Mesh Independence of a Majorized ABCD Method for Sparse PDE-constrained Optimization Problems\thanks{The research of Defeng Sun was supported in part  by the Hong Kong Research Grant Council grant
PolyU153014/18p. The research of Kim-Chuan Toh was supported in part by
the Academic Research Fund (Grant R-146-000-257-112) of
the Ministry of
Education, Singapore.}
}

\titlerunning{Mesh Independence of mABCD for Sparse PDE-constrained Optimization}        

\author{Xiaoliang Song         \and
        Defeng Sun \and
        Kim-Chuan Toh
}


\institute{Xiaoliang Song \at
              Department of Applied Mathematics, The Hong Kong Polytechnic University, Hung Hom, Hong Kong \\
              \email{xiaoliang.song@polyu.edu.hk}           
           \and
           Defeng Sun \at
              Department of Applied Mathematics, The Hong Kong Polytechnic University, Hung Hom, Hong Kong\\
              \email{defeng.sun@polyu.edu.hk} \\
           \and
             Kim-Chuan Toh \at
              Department of Mathematics, and Institute of Operations Research and Analytics, National University of Singapore, 10 Lower Kent Ridge Road, Singapore 119076 \\
              \email{mattohkc@nus.edu.sg}
}

\date{Received: date / Accepted: date}

\maketitle

\begin{abstract}
A majorized accelerated block coordinate descent (mABCD) method in Hilbert space is analyzed to solve {a} sparse PDE-constrained optimization problem via its dual. The finite element approximation method is investigated. The attractive $O(1/k^2)$ iteration complexity of {the mABCD} method for the dual objective function values can be achieved. Based on the convergence result, {we prove the robustness with respect to the mesh size $h$ for {the mABCD method by  establishing} that asymptotically the infinite dimensional ABCD method and finite dimensional discretizations have the same convergence property, and the number of iterations of mABCD method remains almost constant as the discretization is refined.}

\keywords{PDE-constrained optimization\and Sparsity\and Duality approach\and Accelerated block coordinate descent\and Iteration complexity\and Mesh independence}
\end{abstract}

\section{Introduction}
\label{intro}
Optimization problems with constraints which require the solution of a
partial differential equation (PDE) arise widely in many areas of the sciences and
engineering, in particular in problems of design. The development, analysis and implementation of efficient and robust numerical techniques for PDE constrained optimization is of utmost importance for the optimal control
of processes and the optimal design of structures and systems in modern technology. In recent years, a high level of sophistication has been reached for PDE constrained optimization. We refer to the contributions in \cite{HiPiUl,WaWa,BeItKu,PreconditioningforL1control} and many further references given therein.

In this paper, we shall focus on the efficient numerical methods to solve the following elliptic PDE-constrained optimization problem with $L^1$-control cost
\begin{equation}\label{eqn:orginal problems}
           \qquad \left\{ \begin{aligned}
        &\min \limits_{(y,u)\in Y\times U}^{}\ \ J(y,u)=\frac{1}{2}\|y-y_d\|_{L^2(\Omega)}^{2}+\frac{\alpha}{2}\|u\|_{L^2(\Omega)}^{2}+\beta\|u\|_{L^1(\Omega)} \\
        &\qquad{\rm s.t.}\qquad Ly=u+y_r\ \ \mathrm{in}\  \Omega, \\
         &\qquad \qquad \qquad y=0\qquad\quad  \mathrm{on}\ \partial\Omega,\\
         &\qquad \qquad\qquad u\in U_{ad}=\{v(x)|a\leq v(x)\leq b,\ {\rm a.e. }\  \mathrm{on}\ \Omega\}\subseteq U,
                          \end{aligned} \right.\tag{$\mathrm{P}$}
 \end{equation}
where $Y:=H_0^1(\Omega)$, $U:=L^2(\Omega)$, $\Omega\subseteq \mathbb{R}^n$ ($n=2$ or $3$) is a convex, open and bounded domain with $C^{1,1}$- or polygonal boundary $\Gamma$; the desired state $y_d\in L^2(\Omega)$ and the source term $y_r \in L^2(\Omega)$ are given; and $a\leq0\leq b$ and $\alpha$, $\beta>0$. Moreover, the operator $L$ is a second-order linear elliptic differential operator.

It is well-known that {adding the $L^1$-norm penalty} can lead to a sparse optimal control, i.e., the optimal control with small support, which is desirable, for instance, in actuator placement problems \cite{Stadler}. In optimal control of distributed parameter systems, it may be difficult or undesirable to place control devices all over the control domain. Instead, we can decide to localize controllers in small and effective regions. Thus, solving the control problem with an $L^1$-norm {penalty on} the control will give us information about the optimal location {to place} the control devices.

{Throughout this paper, the elliptic PDE is given in the following form}
\begin{equation}\label{eqn:state equations}
  \begin{array}{ccccc}
&Ly=u+y_r  &\mathrm{in} && \Omega \\
&y=0  &\mathrm{on} & &\partial\Omega
  \end{array}
\end{equation}
which satisfies the following assumption.
\begin{assumption}\label{equ:assumption:1}
The linear second-order differential operator $L$ is defined by
 \begin{equation}\label{operator A}
   (Ly)(x):=-\sum \limits^{n}_{i,j=1}\partial_{x_{j}}(a_{ij}(x)y_{x_{i}})+c_0(x)y(x),
 \end{equation}
where functions $a_{ij}(x), c_0(x)\in L^{\infty}(\Omega)$, $c_0\geq0$. Moreover, it is uniformly elliptic, i.e. $a_{ij}(x)=a_{ji}(x)$ and there is a constant $\theta>0$ such that
\begin{equation}\label{equ:operator A coercivity}
  \sum \limits^{n}_{i,j=1}a_{ij}(x)\xi_i\xi_j\geq\theta\|\xi\|^2, \qquad \mathrm{for\ a.a.}\ x\in \Omega\  \mathrm{and}\  \forall \xi \in \mathbb{R}^n.
\end{equation}
{In the above, $y_{x_i}$ denotes the partial derivative of $y(\cdot)$ with respect to $x_i$.}
\end{assumption}

The weak formulation of (\ref{eqn:state equations}) is given by
\begin{equation}\label{eqn:weak form}
  \mathrm{Find}\ y\in H_0^1(\Omega):\ a(y,v)={\langle u+y_r,v\rangle_{L^2(\Omega)}},\quad \forall v \in H_0^1(\Omega),
\end{equation}
with the bilinear form
\begin{equation}\label{eqn:bilinear form}
  a(y,v)=\int_{\Omega}(\sum \limits^{n}_{i,j=1}a_{ji}y_{x_{i}}v_{x_{j}}+c_0yv)\mathrm{d}x,
\end{equation}
or in short $ Ay=B(u+y_r)$,
where $A\in \mathcal{L}(Y,Y^*)$ is the operator induced by the bilinear form $a$, i.e., $Ay=a(y,\cdot)$ and $B\in \mathcal{L}(U,Y^*)$ is defined by $Bu={\langle u,\cdot\rangle_{L^2(\Omega)}}$. Since the bilinear form $a(\cdot,\cdot)$ is symmetric and $U,Y$ are Hilbert spaces, we have $A^*\in\mathcal{L}(Y,Y^*)=A$, and $B^*\in\mathcal{L}(Y,U)$ with $B^*v=v, \forall v\in Y$. 

\begin{remark}\label{more general case}
Although we assume that the Dirichlet boundary condition $y=0$ holds, it should be noted that the assumption is not a restriction and our considerations can also carry over to the more general boundary conditions of Robin type:
\begin{equation*}
  \frac{\partial y}{\partial \nu}+\gamma y=g \quad {\rm on}\  \partial\Omega,
\end{equation*}
where $g\in L^2(\partial\Omega)$ is given and $\gamma\in L^{\infty}(\partial\Omega)$ is a nonnegative coefficient.
\end{remark}

Let us mention some existing numerical methods for solving problem (\ref{eqn:orginal problems}). {For the nonsmooth  problem (\ref{eqn:orginal problems}), semismooth Newton (SSN) methods are  {the primary} choices  in consideration of their locally superlinear convergence, see \cite{Ulbrich1,Ulbrich2,HiPiUl} for more details.} With no doubt, employing the SSN method can derive the solution with {a high accuracy}. However, it should be mentioned that the total error of numerically solving the PDE-constrained problem contains two parts: the discretization error and the iteration error resulted from an algorithm of solving the discretized problem. Obviously, the discretization error accounts for the main part of the total error due to the error order of $O(h)$. Thus, {with the precision of discretization error in mind, algorithms {for very accurately solving the discretized problem} may not  reduce the order of the total error but may incur extra computations.}

{As one may have observed, for finite dimensional large scale optimization problems, some efficient first-order algorithms, such as iterative soft thresholding algorithms (ISTA), accelerated proximal gradient (APG)-based method, alternating direction method of multipliers (ADMM), etc, have become very popular in situations when high accuracy is not sought, see \cite{Blumen,inexactAPG,Beck,Toh,Fazel,SunToh1,SunToh2} and {the} references therein.} 
Hence, employing fast and efficient first-order algorithms with the aim of solving problem (\ref{eqn:orginal problems}) to moderate accuracy is a wise choice. Motivated by the success of some first-order optimization algorithms for finite dimensional optimization problems, to solve problem (\ref{eqn:orginal problems}), the authors \cite{iwADMM} employ an inexact semi-proximal ADMM (isPADMM)
algorithm designed in \cite{SunToh1}. Recently, an APG method was proposed to solve (\ref{eqn:orginal problems}) in \cite{FIP},
{which has the highly desirable iteration complexity of $O(1/k^2)$.}

As far as we know, most of the aforementioned papers are devoted to solve the primal problem.
However, when the primal problem (\ref{eqn:orginal problems}) is discretized by the piecewise linear finite elements and directly solved by some algorithms mentioned above, e.g., SSN, isPADMM and APG, the resulting discretized $L^1$-norm
\begin{eqnarray*}\label{equ:discrete norm}
\|u_h\|_{L^1(\Omega_h)}&=&\int_{\Omega_h}\big{|}\sum\limits_{i=1}^{N_h}u_i\phi_i(x)\big{|}\mathrm{d}x,
\end{eqnarray*}
does not have a decoupled form. To overcome the difficulties, one approach in \cite{WaWa,iwADMM} is introduced by employing an alternative discretization of the $L^1$-norm
\begin{eqnarray*}
\|u_h\|_{L^{1}_h(\Omega_h)}&:=&\sum_{i=1}^{n}|u_i|\int_{\Omega_h}\phi_i(x)\mathrm{d}x.
\end{eqnarray*}
{For} the approximate $L^1$-norm, the authors proved that this approximation technique will not {change} the order of finite element error estimates. Another approach is introduced by Song, Chen and Yu in \cite{mABCDSOPT} by proposing a duality-based approach for solving the problem (\ref{eqn:orginal problems}). Taking advantage of the structure of the dual problem, the authors proposed an inexact symmetric Gauss-Seidel based majorized ABCD (sGS-imABCD) method to solve the discretized dual problem. It should be emphasized that the design of this method combines an inexact 2-block majorized accelerated block coordinate descent (mABCD) method proposed by Cui in \cite{CuiYing} and the recent advances in the inexact symmetric Gauss-Seidel (sGS) decomposition technique developed in \cite{SunToh2,SunToh3}.

In this paper, we will continue to focus on the majorized ABCD algorithm. As known to us, the majorized ABCD method was originally developed for finite dimensional problems. However, when the majorized ABCD algorithm {is}
applied to optimization problems with PDE constraints, some new aspects become important. In particular, a key issue should be considered is how {various} measures of the convergence behavior of the iteration sequence vary with the level of approximation. Such questions come under the category of mesh-independence results. It should be pointed out that mesh independence allows us to predict the convergence of the method {when} applied to the discretized problem {after it} has been analyzed for the infinite dimensional problem. Further, it can be used to improve the performance of the method. Specifically, we can use a prolongated solution on a coarse grid as a good initialization for a finer discretization, which leads to mesh-refinement strategies. Mesh-independence is a theoretical justification for mesh-refinement strategies. More importantly, in \cite{mABCDSOPT}, the numerical results in terms of {the} iteration numbers of mABCD method show that the majorized ABCD method is robust with respect to the mesh size $h$. This phenomenon gives us the strong motivation to establish the mesh independence of the majorized ABCD method, which is the main contribution of this paper.

{ To achieve our goal, 
we first apply the majorized ABCD algorithm on the continuous level for solving the infinite dimensional dual problem of (\ref{eqn:orginal problems}).} Specifically, we will first give a framework of the majorized ABCD algorithm in function space to focus the presentation on structural aspects inherent in the majorized ABCD algorithm and analyze its convergence property. {Then, for the purpose of numerical implementation, a finite element discretized version of the majorized ABCD algorithm is proposed.} Finally, {comparing} the convergence results of the the majorized ABCD algorithm in function space and the discretized version of the majorized ABCD algorithm, {one type of mesh independence for the majorized ABCD method is given.} The result shows that the iteration number $k$ after which the difference $\Phi_h(z^k)-\inf\Phi_h(z)$ has been identified up to less than $\epsilon$ is independent of the mesh size $h$. In other words, we will show that the ``discretized" convergence factor $\tau_h$ defined in the convergence theorem can be bounded by the ``continuous" convergence factor $\tau$.

The remainder of the paper is organized as follows. In Section \ref{sec:2}, we give a majorized accelerate block coordinate descent (mABCD) method in Hilbert space. For the purpose of numerical implementation, in Section \ref{sec:3} the finite element approximation is introduced and the finite element discretizations of the mABCD method is also given.
In Section \ref{sec:4}, we show the mesh independence result of the mABCD method for the sparse PDE-constrained optimization problem (\ref{eqn:orginal problems}). Finally, we conclude our paper in Section \ref{sec:5}.

\section{Duality-based approach}
\label{sec:2}

In this section, we will introduce the duality-based approach to solve problem (\ref{eqn:orginal problems}). First, we will give the dual problem of (\ref{eqn:orginal problems}). Then, to solve the dual problem, we will propose a framework of the majorized ABCD algorithm in function space and focus the presentation on the structural aspects inherent in the majorized ABCD algorithm.

\subsection{\textbf{Dual of problem (\ref{eqn:orginal problems})}}\label{sec:2.1}

{With simple calculations, the dual of problem (\ref{eqn:orginal problems}) can be written}, in its equivalent minimization form, as
\begin{equation}\label{eqn:dual problem}
\begin{aligned}
\min\ \Phi(\lambda,p,\mu):=&{\frac{1}{2}\|A^*p-y_d\|_{L^2(\Omega)}^2}+ \frac{1}{2\alpha}\|p-\lambda-\mu\|_{L^2(\Omega)}^2+\langle p, y_r\rangle_{L^2(\Omega)}\\
&+\delta_{\beta B_{\infty}(0)}(\lambda)+\delta^*_{U_{ad}}(\mu)-\frac{1}{2}\|y_d\|_{L^2(\Omega)}^2,\end{aligned}\tag{$\mathrm{D}$}
\end{equation}
where $p\in H^1_0(\Omega)$, $\lambda,\mu\in L^2(\Omega)$, $B_{\infty}(0):=\{\lambda\in L^2(\Omega): \|\lambda\|_{L^\infty(\Omega)}\leq 1\}$, and for any given nonempty, closed convex subset $C$ of $L^2(\Omega)$, $\delta_{C}(\cdot)$ is the indicator function of $C$. Based on the $L^2$-inner product, we define the conjugate of $\delta_{C}(\cdot)$ as follows:
\begin{equation*}
 \delta^*_{C}(s^*)=\sup\limits_{s\in C}^{}{\langle s^*,s\rangle}_{L^2(\Omega)}.
\end{equation*}
Obviously, by choosing $v=(\lambda, p)$, $w=\mu$ and taking
\begin{eqnarray}
  f(v) &=&\delta_{\beta B_{\infty}(0)}(\lambda)+{\frac{1}{2}\|A^*p-y_d\|_{L^2(\Omega)}^2}+\langle p, y_r\rangle_{L^2(\Omega)}-\frac{1}{2}\|y_d\|_{L^2(\Omega)}^2\label{f function for D},\\
  g(w) &=& \delta^*_{U_{ad}}(\mu)\label{g function for D}, \\
  \phi(v, w) &=&\frac{1}{2\alpha}\|p-\lambda-\mu\|_{L^2(\Omega)}^2\label{phi function for D},
\end{eqnarray}
it is quite clear that our dual problem (\ref{eqn:dual problem})  belongs to a general class of unconstrained, multi-block convex optimization problems with coupled objective function, that is
\begin{equation}\label{eqn:model problem}
\begin{aligned}
\min_{v, w} \theta(v,w):= f(v)+ g(w)+ \phi(v, w),
\end{aligned}
\end{equation}
where $f: \mathcal{V}\rightarrow (-\infty, +\infty ]$ and $g: \mathcal{W}\rightarrow  (-\infty, +\infty ]$ are two convex functions (possibly nonsmooth), $\phi: \mathcal{V}\times \mathcal{W}\rightarrow  (-\infty, +\infty ]$ is a smooth convex function, and $\mathcal{V}$, $\mathcal{W}$ are real Hilbert spaces.
Thus taking advantage of the structure of the dual problem, we will aim to present an algorithm {to solve} problem (\ref{eqn:dual problem}) efficiently.

\subsection{\bf A majorized ABCD algorithm for the general problem (\ref{eqn:model problem})}\label{sec:2.2}

Thanks to the structure of (\ref{eqn:model problem}), Cui in \cite{CuiYing} proposed a majorized accelerate block coordinate descent (mABCD) method. {We give} a brief sketch of mABCD method below. To deal with the general model (\ref{eqn:model problem}), we need some more assumptions on $\phi$.

\begin{assumption}\label{assumption on differentiable with Lipschitz continuous gradients}
The convex function $\phi: \mathcal{V}\times \mathcal{W}\rightarrow  (-\infty, +\infty ]$ is continuously differentiable with Lipschitz continuous gradients.
\end{assumption}

Let us denote $z:=(v, w)\in \mathcal{V}\times\mathcal{W}$. The authors \cite[Theorem 2.3]{Hiriart1984Generalized} provide a second order Mean-Value Theorem for $\phi$, which states that for any $z'$  and $z$ in $\mathcal{V}\times \mathcal{W}$, there exist $ z''\in [z',z]$ and a self-adjoint positive semidefinite operator $\mathcal{G}\in \partial^2\phi(z'')$ such that
\begin{equation*}
\phi(z)= \phi(z')+ \langle\nabla\phi(z'), z- z'\rangle+ \frac{1}{2}\|z'- z\|_{\mathcal{G}}^2,
\end{equation*}
where $\partial^2\phi(z'')$ denotes the Clarke's generalized Hessian at given $z''$ and $[z',z]$ denotes the
the line segment connecting $z'$ and $z$. Under Assumption \ref{assumption on differentiable with Lipschitz continuous gradients}, it is obvious that there exist two self-adjoint positive semidefinite linear operators $\mathcal{Q}$ and $\widehat{\mathcal{Q}}: \mathcal{V}\times\mathcal{W}\rightarrow \mathcal{V}\times\mathcal{W}$ such that for any
$z\in \mathcal{V}\times\mathcal{W}$,
$\mathcal{Q}\preceq\mathcal{G}\preceq\widehat{\mathcal{Q}}$. Thus, for any $z, z'\in \mathcal{V}\times\mathcal{W}$, it holds that
\begin{equation*}
\phi(z)\geq \phi(z')+ \langle\nabla\phi(z'), z- z'\rangle+ \frac{1}{2}\|z'- z\|_{\mathcal{Q}}^2,
\end{equation*}
and
\begin{equation*}
\phi(z)\leq \hat{\phi}(z; z'):= \phi(z ')+ \langle\nabla\phi(z '), z- z '\rangle+ \frac{1}{2}\|z'- z\|_{\widehat{\mathcal{Q}}}^2.
\end{equation*}
Furthermore, we decompose the operators $\mathcal{Q}$ and $\widehat{\mathcal{Q}}$ into the following block structures:
\begin{equation*}
\mathcal{Q}z:=\left(
\begin{array}{cc}
\mathcal{Q}_{11} & \mathcal{Q}_{12}\\
\mathcal{Q}_{12}^* &\mathcal{Q}_{22}
\end{array}
\right)
\left(
\begin{array}{c}
v\\
w
\end{array}
\right),\quad
\widehat{\mathcal{Q}}z:=\left(
\begin{array}{cc}
\widehat{\mathcal{Q}}_{11} & \widehat{\mathcal{Q}}_{12}\\
\widehat{\mathcal{Q}}_{12}^* & \widehat{\mathcal{Q}}_{22}
\end{array}
\right)
\left(
\begin{array}{c}
v\\
w
\end{array}
\right),
\quad\forall z=(v, w)\in \mathcal{U}\times \mathcal{V},
\end{equation*}
and assume $\mathcal{Q}$, $\widehat{\mathcal{Q}}$ satisfy the following assumption.
\begin{assumption}{\rm\textbf{\cite[Assumption 3.1]{CuiYing}}}\label{assumption majorized}
There exist two self-adjoint positive semidefinite linear operators $\mathcal{D}_1: \mathcal{U}\rightarrow \mathcal{U}$ and $\mathcal{D}_2: \mathcal{V}\rightarrow \mathcal{V}$ such that
\begin{equation*}
\widehat{\mathcal{Q}}:=\mathcal{Q}+ {\rm Diag}(\mathcal{D}_1,\mathcal{D}_2).
\end{equation*}
Furthermore, $\widehat{\mathcal{Q}}$ satisfies that $\widehat{\mathcal{Q}}_{11}\succ 0$ and $\widehat{\mathcal{Q}}_{22}\succ 0$.
\end{assumption}

\begin{remark}\label{choice of Hessian}
It is important to note that Assumption {\rm\ref{assumption majorized}} is a realistic assumption in practice. For example, when $\phi$ is a quadratic function, we could choose $\mathcal{Q}=\mathcal{G}=\nabla^2\phi$. If we have $\mathcal{Q}_{11}\succ0$ and $\mathcal{Q}_{22}\succ0$, then Assumption {\rm\ref{assumption majorized}} holds automatically. We should point out that $\phi$ is a quadratic function for many problems in {practical applications}. Fortunately, it should be noted that the function $\phi$ defined in {\rm(\ref{phi function for D})} for our problem {\rm(\ref{eqn:dual problem})} is quadratic and thus we can choose $\mathcal{Q}=\nabla^2\phi$.
\end{remark}

We can now present the majorized ABCD algorithm for (\ref{eqn:model problem}) as follows.

\begin{algorithm}[H]\caption{\textbf{(A majorized ABCD algorithm for (\ref{eqn:model problem}))}}\label{algo1:imabcd}
    \textbf{Input}:{$(\tilde{v}^1,\tilde{w}^1)=({v}^0, {w}^0)\in \textrm{dom} (f)\times \textrm{dom}(g)$. Set $t_1=1$, $k=1$.}\\
  \textbf{Output}:{$ ({v}^k, {w}^k)$}
\begin{description}
  \item[Step 1] Compute
  \begin{equation*}
  \left\{
  \begin{aligned}
    &{v}^{k}
        =\arg\min_{v \in \mathcal{V}}\{f(v)+ \hat{\phi}(v,\tilde{w}^k; \tilde{z}^k)\},\\
    &{w}^{k}
    	=\arg\min_{w \in \mathcal{W}}\{g(w)+ \hat{\phi}(v^k,w; \tilde{z}^k)\},
  \end{aligned}
  \right.
  \end{equation*}
  where $\tilde{z}^k=(\tilde{v}^k,\tilde{w}^k)$.
  \item[Step 2] Set $t_{k+1}=\frac{1+\sqrt{1+4t_k^2}}{2}$ and $\beta_k=\frac{t_k-1}{t_{k+1}}$, compute
  \begin{equation*}
                \tilde{v}^{k+1}=v^{k}+ \beta_{k}(v^{k}-v^{k-1}), \quad  \tilde{w}^{k+1}= w^{k}+ \beta_{k}(w^{k}-w^{k-1}).
  \end{equation*}
  \item[\bf Step 3] If a termination criterion is not met, set $k:=k+1$ and go to Step 1
\end{description}
\end{algorithm}
Here we state the convergence result. For the detailed proof, one could see \cite{CuiYing}. This theorem builds a solid foundation for our subsequent proposed algorithm.

\begin{theorem}{\rm\textbf{\cite[Theorem 3.2]{CuiYing}}}\label{imABCD convergence}
Suppose that Assumption {\rm\ref{assumption majorized}} holds and the solution set $\Omega$ of the problem {\rm(\ref{eqn:model problem})} is non-empty. Let
$z^*=(v^*,w^*)\in \Omega$. Then the sequence $\{{z}^k\}:=\{({v}^k,{w}^k)\}$ generated by the Algorithm {\rm\ref{algo1:imabcd}} satisfies that
\begin{equation*}
\theta({z}^k)- \theta(z^*)\leq \frac{2\|{z}^0- z^*\|_{\mathcal{S}}^2}{(k+1)^2} \quad \forall k\geq 1,
\end{equation*}
where $\theta(\cdot)$ is the objective function of {\rm(\ref{eqn:model problem})} and $\mathcal{S}:={\rm{Diag}}(\mathcal{D}_1,\mathcal{D}_2+\mathcal{Q}_{22})$.
\end{theorem}

\subsection{\textbf{The sGS-majorized ABCD method in Hilbert Space for (\ref{eqn:dual problem})}}\label{sec:2.3}

Now, we can apply Algorithm \ref{algo1:imabcd} to (\ref{eqn:dual problem}), where $(\lambda,p)$ is taken as one block, and $\mu$ is taken as the other one. Let us denote $z=(\lambda, p,\mu)$. Since $\phi$ defined in (\ref{phi function for D}) for (\ref{eqn:dual problem}) is quadratic, we can take
\begin{equation*}
\mathcal{Q}:=
\frac{1}{\alpha}\left(
\begin{array}{ccc}
\mathcal{I} &\quad -\mathcal{I} &\quad \mathcal{I}\\
-\mathcal{I} & \quad  \mathcal{I}&\quad -\mathcal{I}\\
\mathcal{I} & \quad-\mathcal{I}&\quad\mathcal{I}
\end{array}
\right),
\end{equation*}
where
\begin{equation*}
\mathcal{Q}_{11}:=
\frac{1}{\alpha}\left(
\begin{array}{cc}
\mathcal{I}& \quad-\mathcal{I}\\
-\mathcal{I} & \quad \mathcal{I}
\end{array}
\right), \quad \mathcal{Q}_{22}:= \frac{1}{\alpha}\mathcal{I}.
\end{equation*}
Additionally, we assume that there exist two self-adjoint positive semidefinite operators $\mathcal{D}_1$ and $\mathcal{D}_2$, such that Assumption \ref{assumption majorized} holds. Thus, it implies that we should majorize $\phi(\lambda,p,\mu)$ at $z'=(\lambda', p',\mu')$ as
\begin{equation}\label{majorized function}
\begin{aligned}
 \phi(z) \leq \hat{\phi}(z;z')
:=&\frac{1}{2\alpha}\|-p+\lambda+\mu\|_{L^2(\Omega)}^2+\frac{1}{2} \Big\langle        \left(\begin{array}{c}
          \lambda-\lambda' \\
          p-p'
        \end{array}\right),
        \mathcal{D}_1\left(\begin{array}{c}
          \lambda-\lambda' \\
          p-p'
        \end{array}\right)
        \Big\rangle_{L^2(\Omega)}\\
        &+\frac{1}{2}\big\langle\mu-\mu',\mathcal{D}_2(\mu-\mu')\big\rangle_{L^2(\Omega)}.
 \end{aligned}
\end{equation}
Thus, the framework of mABCD for (\ref{eqn:dual problem}) is given below:

\begin{algorithm}[H]
  \caption{\textbf{(mABCD algorithm for (\ref{eqn:dual problem}))}}\label{algo1:ABCD algorithm for (D)}
 \textbf{Input}:{$(\tilde{\lambda}^1, \tilde{p}^1,\tilde{\mu}^1 )=({\lambda}^0, {p}^0,\mu^0)\in [-\beta,\beta]\times H^{1}_0(\Omega)\times L^2(\Omega)$. $\mathcal{T}\succeq0$. Set $k= 1, t_1= 1.$}

 \textbf{Output}:{$ ({\lambda}^k, {p}^k,{\mu}^k)$}
\begin{description}
  \item[Step 1]
Compute 
\begin{eqnarray*}
 ({\lambda}^{k},p^k)
        &=&\arg\min\delta_{[-\beta,\beta]}(\lambda)+\frac{1}{2}\|A^*p-y_{d}\|_{L^2(\Omega)}^2+\langle p, y_r\rangle_{L^2(\Omega)}+\frac{1}{2\alpha}\|-p+\lambda+\tilde{\mu}^k\|_{L^2(\Omega)}^2\\
        &&\qquad\qquad+\frac{1}{2} \Big\langle        \left(\begin{array}{c}
          \lambda-\tilde{\lambda}^k \\
          p-\tilde{p}^k
        \end{array}\right),
        \mathcal{D}_1\left(\begin{array}{c}
          \lambda-\tilde{\lambda}^k \\
          p-\tilde{p}^k
        \end{array}\right)
        \Big\rangle,\\
  {\mu}^{k}&=&\arg\min\frac{1}{2\alpha}\|\mu-(p^k-\lambda^k)\|_{L^2(\Omega)}^2+\delta^*_{[a,b]}(\mu)+\frac{1}{2}\langle\mu-\tilde{\mu}^k,\mathcal{D}_2(\mu-\tilde{\mu}^k)\rangle.
\end{eqnarray*}
  \item[Step 2] Set $t_{k+1}=\frac{1+\sqrt{1+4t_k^2}}{2}$ and $\beta_k=\frac{t_k-1}{t_{k+1}}$, compute
\begin{eqnarray*}
\tilde{\lambda}^{k+1}= {\lambda}^{k}+ \beta_{k}({\lambda}^{k}-{\lambda}^{k-1}),\quad
 \tilde {p}^{k+1}={p}^{k}+ \beta_{k}({p}^{k}-{p}^{k-1}), \quad
 \tilde {\mu}^{k+1}={\mu}^{k}+ \beta_{k}({\mu}^{k}-{\mu}^{k-1}).
\end{eqnarray*}
\item[\bf Step 3] If a termination criterion is not met, set $k:=k+1$ and go to Step 1
\end{description}
\end{algorithm}

We now can discuss the issue on how to choose two operators $\mathcal{D}_1$ and $\mathcal{D}_2$ for Algorithm \ref{algo1:ABCD algorithm for (D)}. As we know, choosing {the operators $\mathcal{D}_1$ and $\mathcal{D}_2$ appropriately is important for numerical computation}. Note that for numerical efficiency, the general principle is that both $\mathcal{D}_1$ and $\mathcal{D}_2$ should be chosen as small as possible such that $({\lambda}^{k},{p}^{k})$ and ${\mu}^{k}$ could take larger step-lengths while the corresponding subproblems still can be solved relatively easily.

Firstly, for the proximal term $\frac{1}{2}\|\mu-\tilde{\mu}^k\|^2_{\mathcal{D}_2}$, since $\mathcal{Q}_{22}=\frac{1}{\alpha}\mathcal{I}\succ0$, we can choose $\mathcal{D}_2=0$. Then, it is obvious that the optimal solution of the $\mu$-subproblem at $k$-th iteration is unique and also has a closed form {solution given by}
\begin{equation}\label{closed form solution of mu}
  \mu^k=p^k-\lambda^k-\alpha{\rm\Pi}_{[a,b]}(\frac{1}{\alpha}(p^k-\lambda^k)).
\end{equation}
Next, we focus on how to choose $\mathcal{D}_1$. Ignoring the proximal term
\begin{equation*}
  \frac{1}{2}\Big\langle\left(\begin{array}{c}
          \lambda-\tilde{\lambda}^k \\
          p-\tilde{p}^k
        \end{array}\right),
        \mathcal{D}_1\left(\begin{array}{c}
          \lambda-\tilde{\lambda}^k \\
          p-\tilde{p}^k
        \end{array}\right)
       \Big \rangle,
\end{equation*}
it is clear that the subproblem with respect to $(\lambda,p)$ at $k$-th iteration can be equivalently rewritten as:
\begin{equation}\label{sGS subproblem for D}
 \min \delta_{[-\beta,\beta]}(\lambda)+\frac{1}{2}\Big\langle \left(\begin{array}{c}
                                                             \lambda \\
                                                             p
                                                           \end{array}\right)
  ,\mathcal{H}\left(\begin{array}{c}
                                                             \lambda \\
                                                             p
                                                           \end{array}\right)\Big\rangle-
                                                           \Big\langle r,\left(\begin{array}{c}
                                                             \lambda \\
                                                             p
                                                           \end{array}\right) \Big\rangle,
\end{equation}
where $\mathcal{H}=
\left(
\begin{array}{cc}
                        \frac{1}{\alpha}\mathcal{I} & \quad -\frac{1}{\alpha}\mathcal{I} \\
                        -\frac{1}{\alpha}\mathcal{I} & \quad AA^*+\frac{1}{\alpha}\mathcal{I} \\
                      \end{array}
\right)$ and $r=\left(\begin{array}{c}
                                                             -\frac{1}{\alpha}\tilde{\mu}^k \\
                                                             -y_r+Ay_d+\frac{1}{\alpha}\tilde{\mu}^k
                                                           \end{array}\right)$,
whose objective function of (\ref{sGS subproblem for D}) is the sum of a two-block quadratic function and a non-smooth function involving only the first block, thus the symmetric Gauss-Seidel (sGS) technique proposed recently by Li, Sun and Toh \cite{SunToh2,SunToh3}, could be used to solve it. For later discussions, we consider a splitting of any given self-adjoint positive semidefinite linear operator $\mathcal{Q}$
\begin{equation}\label{equ:splitting}
  \mathcal{Q}=\mathcal{D}+\mathcal{U}+\mathcal{U}^*,
\end{equation}
where $\mathcal{U}$ denotes the strict upper triangular part of $\mathcal{Q}$ and $\mathcal{D}$ is the diagonal of $\mathcal{Q}$. Moreover, we assume that $\mathcal{D}\succ0$ and define the following self-adjoint positive semidefinite linear operator
\begin{equation}\label{SGSoperator}
  {\rm sGS}(\mathcal{Q}):=\mathcal{U}\mathcal{D}^{-1}\mathcal{U}^*.
\end{equation}
Thus, to achieve our goal, we choose
\begin{eqnarray*}
  \mathcal{D}_1:&=&\mathrm{sGS}\left(\mathcal{H}\right)
                      =\left(
                      \begin{array}{cc}
                        \frac{1}{\alpha}(\alpha A A^*+\mathcal{I})^{-1} & \quad 0 \\
                        0 & \quad 0 \\
                      \end{array}
                  \right).
\end{eqnarray*}
Then according to \cite[Theorem 2.1]{SunToh3}, solving the $(\lambda,p)$-subproblem
\begin{equation*}
\begin{aligned}
{(\lambda^k,p^k) =}  \mbox{argmin}_{\lambda,p}\ &\delta_{[-\beta,\beta]}(\lambda)+\frac{1}{2}\|A^*p-y_{d}\|_{L^2(\Omega)}^2+\langle p, y_r\rangle_{L^2(\Omega)}+\frac{1}{2\alpha}\|-p+\lambda+\tilde{\mu}^k\|_{L^2(\Omega)}^2\\
       &+\frac{1}{2}\Big\langle        \left(\begin{array}{c}
          \lambda-\tilde{\lambda}^k \\
          p-\tilde{p}^k
        \end{array}\right),
        \mathcal{D}_1\left(\begin{array}{c}
          \lambda-\tilde{\lambda}^k \\
          p-\tilde{p}^k
        \end{array}\right)
        \Big\rangle,
\end{aligned}
\end{equation*}
is equivalent to computing {$(\lambda^k,p^k)$} via the following procedure:
\begin{equation*}
\left\{\begin{aligned}
 &\hat{p}^{k}=\arg\min\frac{1}{2}\|A^* p-y_{d}\|_{L^2(\Omega)}^2+ \frac{1}{2\alpha}\|p-\tilde{\lambda}^k-\tilde{\mu}^k\|_{L^2(\Omega)}^2+\langle p, y_r\rangle_{L^2(\Omega)},\\
 &{\lambda}^{k}
        =\arg\min\frac{1}{2\alpha}\|\lambda-(\hat{p}^{k}-\tilde{\mu}^k)\|_{L^2(\Omega)}^2+\delta_{[-\beta,\beta]}(\lambda),\\
&{p}^{k}=\arg\min\frac{1}{2}\|A^*p-y_{d}\|_{L^2(\Omega)}^2+ \frac{1}{2\alpha}\|p-{\lambda}^k-\tilde{\mu}^k\|_{L^2(\Omega)}^2+\langle p, y_r\rangle_{L^2(\Omega)}.
\end{aligned}\right.
\end{equation*}
\begin{remark}\label{solutions of subproblem}
Specifically, for the $\lambda$-subproblem of Algorithm {\rm\ref{algo1:ABCD algorithm for (D)}} at {the}
$k$-th iteration, it has a closed form solution which {is given} by
\begin{equation*}
\lambda^k={\rm\Pi}_{[-\beta,\beta]}(\hat{p}^{k}-\tilde{\mu}^k).
\end{equation*}
For the $\hat{p}$-subproblem, it is obvious that solving the subproblem is equivalent to solving the following system:
\begin{equation*}
  A(A^*p-y_d)+\frac{1}{\alpha}(p-\tilde{\lambda}^k-\tilde{\mu}^k)+y_r=0.
\end{equation*}
Moreover, to solve the ${p^k}$-subproblem, we only need to replace $\tilde{\lambda}^k$ by ${\lambda}^k$ in the right-hand term. Thus, all the numerical techniques for the block $\hat{p}^k$ is also applicable for the block ${p^k}$. 
\end{remark}

At last, combining a 2-block majorized ABCD and the recent advances in the symmetric Gauss-Seidel (sGS) decomposition technique,
{a sGS} based majorized ABCD (sGS-mABCD) algorithm for (\ref{eqn:dual problem}) is presented as follows.

\begin{algorithm}[H]
  \caption{\textbf{(sGS-mABCD algorithm for (\ref{eqn:dual problem}))}}\label{algo1:sGS-mABCD algorithm for (D)}
 \textbf{Input}:{$(\tilde{\lambda}^1, \tilde{p}^1,\tilde{\mu}^1)=( {\lambda}^0, {p}^0, \mu^0)\in [-\beta,\beta]\times H^{1}_0(\Omega)\times L^2(\Omega)$. Set $k= 1, t_1= 1.$}

 \textbf{Output}:{$ ({\lambda}^k, {p}^k, {\mu}^k)$}

\begin{description}
  \item[Step 1]
Compute 
\begin{eqnarray*}
 \hat{p}^{k}&=&\arg\min\frac{1}{2}\|A^* p-y_{d}\|_{L^2(\Omega)}^2+ \frac{1}{2\alpha}\|p-\tilde{\lambda}^k-\tilde{\mu}^k\|_{L^2(\Omega)}^2+\langle p, y_r\rangle_{L^2(\Omega)},\\
 \\
 {\lambda}^{k}
        &=&\arg\min\frac{1}{2\alpha}\|\lambda-(\hat{p}^{k}-\tilde{\mu}^k)\|_{L^2(\Omega)}^2+\delta_{[-\beta,\beta]}(\lambda),\\
 \\
{p}^{k}&=&\arg\min\frac{1}{2}\|A^*p-y_{d}\|_{L^2(\Omega)}^2+ \frac{1}{2\alpha}\|p-{\lambda}^k-\tilde{\mu}^k\|_{L^2(\Omega)}^2+\langle p, y_r\rangle_{L^2(\Omega)},\\
{\mu}^{k}&=&\arg\min\frac{1}{2\alpha}\|\mu-({p}^k-{\lambda}^k)\|_{L^2(\Omega)}^2+\delta^*_{[a,b]}(\mu).
\end{eqnarray*}

  \item[Step 2] Set $t_{k+1}=\frac{1+\sqrt{1+4t_k^2}}{2}$ and $\beta_k=\frac{t_k-1}{t_{k+1}}$, compute
\begin{eqnarray*}
\tilde{\lambda}^{k+1}= {\lambda}^{k}+ \beta_{k}({\lambda}^{k}-{\lambda}^{k-1}),\quad
 \tilde {p}^{k+1}={p}^{k}+ \beta_{k}({p}^{k}-{p}^{k-1}), \quad
 \tilde {\mu}^{k+1}={\mu}^{k}+ \beta_{k}({\mu}^{k}-{\mu}^{k-1}).
\end{eqnarray*}
\item[\bf Step 3] If a termination criterion is not met, set $k:=k+1$ and go to Step 1
\end{description}
\end{algorithm}

Employing Theorem \ref{imABCD convergence}, we have the following convergence result for Algorithm \ref{algo1:sGS-mABCD algorithm for (D)}.
\begin{theorem}\label{imABCD convergence for dual problem}
Suppose that the solution set $\Theta$ of Problem {\rm(\ref{eqn:dual problem})} is non-empty. Let
$(\lambda^*,p^*,\mu^*)\in \Theta$. Then the sequence $\{(\lambda^k,p^k,\mu^k)\}$ generated by Algorithm {\rm\ref{algo1:sGS-mABCD algorithm for (D)}} satisfies that
\begin{equation}\label{iteration complexity of dual problem for D}
{\Phi}(\lambda^k,p^k,\mu^k)- {\Phi}(\lambda^*,p^*,\mu^*)\leq \frac{4\tau}{(k+1)^2},
\end{equation}
where $\Phi(\cdot)$ is the objective function of the dual problem {\rm(\ref{eqn:dual problem})} and

\begin{eqnarray*}
 &&\tau=\frac{1}{2}\langle\left(\begin{array}{c}
          \lambda^*-\lambda^0 \\
          p^*-p^0\\
          \mu^*-\mu^0 \\
        \end{array}\right),
        \mathcal{S}\left(\begin{array}{c}
          \lambda^*-\lambda^0 \\
          p^*-p^0\\
          \mu^*-\mu^0 \\
        \end{array}\right)
        \rangle,\quad \mathcal{S}=
\left(
  \begin{array}{ccc}
    \frac{1}{\alpha}(\alpha A^*A+\mathcal{I})^{-1} & 0 & \quad0 \\
    0 & 0 & \quad0 \\
    0 & 0 & \quad\frac{1}{\alpha}\mathcal{I} \\
  \end{array}
\right).
\end{eqnarray*}
\end{theorem}

\section{Finite element discretization}
\label{sec:3}

\subsection{\textbf {Piecewise linear finite elements discretization}}
\label{sec:3.1}

To numerically solve problem (\ref{eqn:orginal problems}), we consider the finite element method, in which the state $y$ and the control $u$ are both discretized by the piecewise linear, globally continuous finite elements. To achieve this aim, let us fix the assumptions on the discretization by finite elements. We first consider a family of regular and quasi-uniform triangulations $\{\mathcal{T}_h\}_{h>0}$ of $\bar{\Omega}$. For each cell $T\in \mathcal{T}_h$, let us define the diameter of the set $T$ by $\rho_{T}:={\rm diam}\ T$ and define $\sigma_{T}$ to be the diameter of the largest ball contained in $T$. The mesh size of the grid is defined by $h=\max_{T\in \mathcal{T}_h}\rho_{T}$. We suppose that the following regularity assumption on the triangulation {is} satisfied, which {is} standard in the context of error estimates.

\begin{assumption}[regular and quasi-uniform triangulations]\label{regular and quasi-uniform triangulations}
There exist two positive constants $\kappa$ and $\tau$ such that
   \begin{equation*}
   \frac{\rho_{T}}{\sigma_{T}}\leq \kappa,\quad \frac{h}{\rho_{T}}\leq \tau,
 \end{equation*}
hold for all $T\in \mathcal{T}_h$ and all $h>0$. Moreover, let us define $\bar{\Omega}_h=\bigcup_{T\in \mathcal{T}_h}T$, and let ${\Omega}_h \subset\Omega$ and $\Gamma_h$ denote its interior and its boundary, respectively. In the case that $\Omega$ is a convex polyhedral domain, we have $\Omega=\Omega_h$. In the case {that} $\Omega$ has a $C^{1,1}$- boundary $\Gamma$, we assume that $\bar{\Omega}_h$ is convex and all boundary vertices of $\bar{\Omega}_h$ are contained in $\Gamma$, such that
\begin{equation*}
  |\Omega\backslash {\Omega}_h|\leq c h^2,
\end{equation*}
where $|\cdot|$ denotes the measure of the set and $c>0$ is a constant.
\end{assumption}

On account of the homogeneous boundary condition of the state equation, we use
\begin{equation}\label{state discretized space}
  Y_h =\left\{y_h\in C(\bar{\Omega})~\big{|}~y_{h|T}\in \mathcal{P}_1~ {\rm{for\ all}}~ T\in \mathcal{T}_h~ \mathrm{and}~ y_h=0~ \mathrm{in } ~\bar{\Omega}\backslash {\Omega}_h\right\}
\end{equation}
as the discrete state space, where $\mathcal{P}_1$ denotes the space of polynomials of degree less than or equal to $1$. As mentioned above, we also use the same discrete space to discretize {the} control $u$, thus we define
\begin{equation}\label{control discretized space}
   U_h =\left\{u_h\in C(\bar{\Omega})~\big{|}~u_{h|T}\in \mathcal{P}_1~ {\rm{for\ all}}~ T\in \mathcal{T}_h~ \mathrm{and}~ u_h=0~ \mathrm{in } ~\bar{\Omega}\backslash{\Omega}_h\right\}.
\end{equation}
For a given regular and quasi-uniform triangulation $\mathcal{T}_h$ with nodes $\{x_i\}_{i=1}^{N_h}$, let $\{\phi_i(x)\} _{i=1}^{N_h}$ be a set of nodal basis functions, which span $Y_h$ as well as $U_h$ and satisfy the following properties:
\begin{eqnarray}\label{basic functions properties}
  &&\phi_i(x) \geq 0, \quad
  \|\phi_i\|_{\infty} = 1 \quad \forall i=1,2,...,N_h,
 \quad \sum\limits_{i=1}^{N_h}\phi_i(x)=1, \forall x\in \Omega_h.
\end{eqnarray}
The elements $u_h\in U_h$ and $y_h\in Y_h$ can be represented in the following forms respectively,
\begin{equation*}
  u_h=\sum \limits_{i=1}^{N_h}u_i\phi_i(x),\quad y_h=\sum \limits_{i=1}^{N_h}y_i\phi_i(x),
\end{equation*}
where $u_h(x_i)=u_i$ and $y_h(x_i)=y_i$. Let $U_{ad,h}$ denote the discrete feasible set, which is defined by
\begin{eqnarray*}
  U_{ad,h}:&=&U_h\cap U_{ad}\\
           &=&\left\{z_h=\sum \limits_{i=1}^{N_h} z_i\phi_i(x)~\big{|}~a\leq z_i\leq b, \forall i=1,...,N_h\right\}\subset U_{ad}.
\end{eqnarray*}

From the perspective of numerical implementation, we introduce the following stiffness and mass matrices:
\begin{equation*}
 K_h = \left(a(\phi_i, \phi_j)\right)_{i,j=1}^{N_h},\quad  M_h=\left(\int_{\Omega_h}\phi_i(x)\phi_j(x){\mathrm{d}}x\right)_{i,j=1}^{N_h},
\end{equation*}
and let $y_{r,h}$, $y_{d,h}$ be the projections of $y_r$ and $y_d$ onto $Y_h$, respectively,
\begin{equation*}
  y_{r,h}=\sum\limits_{i=1}^{N_h}y_r^i\phi_i(x),\quad y_{d,h}=\sum\limits_{i=1}^{N_h}y_d^i\phi_i(x).
\end{equation*}

Moreover, for the requirement of the subsequent discretized algorithms, {next we} introduce the lumped mass matrix $W_h$
\begin{equation*}
  W_h:={\rm{diag}}\left(\int_{\Omega_h}\phi_i(x)\mathrm{d}x\right)_{i=1}^{N_h},
\end{equation*}
which is a diagonal matrix, and define an alternative discretization of the $L^1$-norm:
\begin{equation}\label{equ:approximal L1}
  \|u_h\|_{L^{1}_h(\Omega)}:=\sum_{i=1}^{N_h}|u_i|\int_{\Omega_h}\phi_i(x)\mathrm{d}x=\|W_h{\bf u}\|_1,
\end{equation}
which is a weighted $l^1$-norm of the coefficients of $u_h$.
More importantly, the following results about the mass matrix $M_h$ and the lumped mass matrix $W_h$ hold.
\begin{proposition}{\rm{\textbf{\cite[Table 1]{Wathen}}}}\label{eqn:martix properties1}
$\forall$ ${\bf z}=(z_1,z_2,...,z_{N_h})\in \mathbb{R}^{N_h}$, the following inequalities hold:
\begin{eqnarray}
 \label{Winequality1}&&\|{\bf z}\|^2_{M_h}\leq\|{\bf z}\|^2_{W_h}\leq \gamma\|{\bf z}\|^2_{M_h} \quad where \ \gamma=
 \left\{ \begin{aligned}
         &4  \quad if \ n=2, \\
         &5  \quad if \ n=3,
                           \end{aligned} \right.
                           \\
 \label{Winequality2}&&\int_{\Omega_h}|\sum_{i=1}^n{z_i\phi_i(x)}|~\mathrm{d}x\leq\|W_h{\bf z}\|_1.
\end{eqnarray}
\end{proposition}

To analyze the error between $\|u_h\|_{L^{1}_h(\Omega)}$ and $\|u_h\|_{L^{1}(\Omega)}$, we first introduce the nodal interpolation operator $I_h$. For a given regular and quasi-uniform triangulation $\mathcal{T}_h$ of $\Omega$ with nodes $\{x_i\}_{i=1}^{N_h}$, we define
\begin{equation}\label{nodal interpolation operator}
(I_hw)(x)=\sum_{i=1}^{N_h}w(x_i)\phi_i(x) \ {\rm\ for\ any}\ w\in L^1(\Omega).
\end{equation}
{Concerning} the interpolation error estimate, we have the following result, see  {\rm\cite[Theorem 3.1.6]{Ciarlet}} for more details.

\begin{lemma}\label{interpolation error estimate}
For all $w\in W^{k+1,p}(\Omega)$, $k\geq 0$, $p,q\in [0,+\infty)$, and $0\leq m\leq k+1$, we have
\begin{equation}\label{interpolation error estimate inequilaty}
  \|w-I_hw\|_{W^{m,q}(\Omega)}\leq c_I h^{k+1-m}\|w\|_{W^{k+1,p}(\Omega)}.
\end{equation}
\end{lemma}
Thus, according to Lemma \ref{interpolation error estimate}, we have the following error estimate results.
\begin{proposition}\label{eqn:martix properties}
$\forall$ ${\bf z}=(z_1,z_2,...,z_{N_h})\in \mathbb{R}^{N_h}$, let $z_h=\sum\limits_{i=1}^{N_h}z_i\phi_i(x)$, then the following inequalities hold
\begin{eqnarray}
 \label{Winequality3}&&0\leq\|z_h\|_{L^{1}_h(\Omega)}-\|z_h\|_{L^{1}(\Omega)}\leq C\,h\,\|z_h\|_{H^1(\Omega)},
\end{eqnarray}
where $C$ is a constant.
\begin{proof}
Obviously, we have
\begin{eqnarray*}
  \|z_h\|_{L^{1}_h(\Omega)}- {\|z_h\|_{L^{1}(\Omega)}}&=&\int_{\Omega_h}\sum_{i=1}^{N_h}|z_i|\phi_i(x)\mathrm{d}x-\int_{\Omega_h}|\sum_{i=1}^n{z_i\phi_i(x)}|~\mathrm{d}x\\
  &=&\int_{\Omega_h} {\big((I_h|z_h|)(x)-|z_h(x)|\big)} \mathrm{d}x.
\end{eqnarray*}
Moreover, due to $z_h\in U_h$, we have $|z_h|\in H^1(\Omega)$. Thus employing Lemma \ref{interpolation error estimate}, we have
\begin{eqnarray*}
\int_{\Omega_h}\big( (I_h|z_h|)(x)-|z_h(x)|\big) \mathrm{d}x&\leq& c_{\Omega}\|I_h|z_h|- |z_h|\|_{L^2(\Omega)}
  =C \,h\, \|z_h\|_{H^1(\Omega)}.
\end{eqnarray*}
{Thus, the proof is completed.}
\end{proof}
\end{proposition}

\subsection{\textbf{A discretized form of sGS-majorized ABCD algorithm for (\ref{eqn:discretized matrix-vector dual problem})}}\label{sec:3.2}
Although an efficient majorized ABCD algorithm in Hilbert space is presented in Section \ref{sec:2}, for the purpose of numerical implementation, we should give the finite element discretizations of the majorized ABCD method. First, employing the piecewise linear, globally continuous finite elements to discretize all the dual variables, then
a type of finite element discretization of (\ref{eqn:dual problem}) is given as follows
\begin{equation}\label{eqn:discretized matrix-vector dual problem}
\begin{aligned}
\min\limits_{{\bm \lambda},{\bf p},{\bm \mu}\in \mathbb{R}^{N_h}}
\Phi_h({\bm \lambda},{\bf p},{\bm \mu}):=&
\frac{1}{2}\|K_h {\bf p}-{M_h}{\bf y_{d}}\|_{M_h^{-1}}^2+ \frac{1}{2\alpha}\|{\bm \lambda}+{\bm \mu}-{\bf p}\|_{M_h}^2+\langle M_h{\bf y_r}, {\bf p}\rangle\\
&+ \delta_{[-\beta,\beta]}({\bm\lambda})+ \delta^*_{[a,b]}({M_h}{\bm\mu})-\frac{1}{2}\|{\bf y_d}\|^2_{M_h}.
\end{aligned}\tag{$\mathrm{D_h}$}
\end{equation}
Obviously, by choosing $v=(\bm \lambda, {\bf p})$, $w=\bm \mu$ and taking
\begin{eqnarray}
  f_h(v) &=&\delta_{[-\beta,\beta]}(\bm\lambda)+\frac{1}{2}\|K_h {\bf p}-{M_h}{\bf y_{d}}\|_{M_h^{-1}}^2+\langle M_h {\bf y_r}, {\bf p}\rangle-\frac{1}{2}\|{\bf y_d}\|^2_{M_h} \label{g function for Dh}, \\
  g_h(w) &=&\delta^*_{[a,b]}({M_h}\bm\mu)\label{f function for Dh},\\
  \phi_h(v, w) &=&\frac{1}{2\alpha}\|\bm\lambda- {\bf p} + \bm\mu\|_{M_h}^2\label{phi function for Dh},
\end{eqnarray}
(\ref{eqn:discretized matrix-vector dual problem}) also belongs to the problem of form (\ref{eqn:model problem}). Thus, Algorithm \ref{algo1:imabcd} also can be applied to (\ref{eqn:discretized matrix-vector dual problem}). Let us denote ${\bf z}=({\bm\lambda},{\bf p},{\bm\mu})$. As shown in Section \ref{sec:3.2}, we should first majorize the coupled function {$\phi_h$} defined in (\ref{phi function for Dh}) for (\ref{eqn:discretized matrix-vector dual problem}). Since $\phi_h$ is quadratic, we can take
\begin{equation}\label{Discretized Hessian Matrix}
\mathcal{Q}_h:=
{
\frac{1}{\alpha}\left(
\begin{array}{ccc}
M_h &\quad -M_h &\quad M_h\\
-M_h & \quad  M_h&\quad -M_h\\
M_h & \quad-M_h&\quad M_h
\end{array}
\right),
}
\end{equation}
where
\begin{equation*}
\mathcal{Q}_h^{11}:=
\frac{1}{\alpha}\left(
\begin{array}{cc}
M_h& \quad-M_h\\
-M_h & \quad M_h
\end{array}
\right),\quad
\mathcal{Q}_h^{22}:= \frac{1}{\alpha} M_h.
\end{equation*}
Moreover, we assume that there exist two self-adjoint positive semidefinite operators $D_{1h}$ and $D_{2h}$, which satisfy Assumption \ref{assumption majorized}. Then, we majorize ${\phi_h}({\bm \lambda},{\bf p},{\bm \mu})$ at $z'=({\bm \lambda'}, {\bf p'},{\bm\mu'})$ as
\begin{equation}\label{majorized function2}
\begin{aligned}
 {\phi_h({\bf z})} \leq {\hat{\phi}_h({\bf z};{\bf z'})}
 =& \frac{1}{2\alpha}\|{\bm\lambda}+ {\bm\mu}-{\bf p}\|_{M_h}^2+\frac{1}{2}\left\|\left(\begin{array}{c}
             {\bm\lambda} \\
             {\bf p}
             \end{array}\right)
-\left(\begin{array}{c}
{\bm \lambda'} \\
{\bf p'}
\end{array}\right)\right\|^2_{D_{1h}}+\frac{1}{2}\|{\bm\mu}-\bm \mu'\|^2_{D_{2h}}.
 \end{aligned}
\end{equation}
Thus, the framework of mABCD for (\ref{eqn:discretized matrix-vector dual problem}) is given as {follows.}

\begin{algorithm}[H]\label{algo1:imABCD algorithm for (Dh)}
  \caption{\textbf{(mABCD algorithm for (\ref{eqn:discretized matrix-vector dual problem}))}}
 \textbf{Input}:{$(\tilde{{\bm\lambda}}^1, \tilde{{\bf p}}^1,\tilde{{\bm\mu}}^1)=({\bm \lambda}^0, {\bf p}^0,\bm\mu^0)\in {\rm dom} (\delta^*_{[a,b]})\times [-\beta,\beta]\times \mathbb{R}^{N_h}$. Set $k= 1, t_1= 1.$}

 \textbf{Output}:{$ ({\bm\lambda}^k, {\bf p}^k,{\bm\mu}^k)$}

\begin{description}
  \item[Step 1] Compute 
\begin{eqnarray*}
({\bm\lambda}^{k},{\bf p}^{k})&=&\arg\min\delta_{[-\beta,\beta]}({\bm\lambda})+\frac{1}{2}\|K_h {\bf p}-{M_h}{\bf y_{d}}\|_{M_h^{-1}}^2+\langle M_h {\bf y_r}, {\bf p}\rangle+\frac{1}{2\alpha}\|{\bm \lambda}-{\bf p}+\tilde{{\bm\mu}}^k\|_{M_h}^2\\
&&\qquad\qquad+\frac{1}{2}\left\|\left(\begin{array}{c}
             {\bm\lambda} \\
             {\bf p}
             \end{array}\right)
-\left(\begin{array}{c}
\tilde{{\bm\lambda}}^k \\
\tilde{{\bf p}}^k
\end{array}\right)\right\|^2_{D_{1h}}.\\
  {\bm\mu}^{k}&=&\arg\min\delta^*_{[a,b]}(M_h{\bm\mu})+\frac{1}{2\alpha}\|{\bm\mu}-({\bf p}^k-{\bm\lambda}^k)\|_{M_h}^2+\frac{1}{2}\|{\bm\mu}-\tilde{{\bm\mu}}^k\|^2_{D_{2h}},
\end{eqnarray*}

  \item[Step 2] Set $t_{k+1}=\frac{1+\sqrt{1+4t_k^2}}{2}$ and $\beta_k=\frac{t_k-1}{t_{k+1}}$, Compute
\begin{eqnarray*}
\tilde{\bm\lambda}^{k+1}= {\bm\lambda}^{k}+ \beta_{k}({\bm\lambda}^{k}-{\bm\lambda}^{k-1}),\quad
\tilde{\bf p}^{k+1}={\bf p}^{k}+ \beta_{k}({\bf p}^{k}-{\bf p}^{k-1}), \quad\tilde{\bm\mu}^{k+1}={\bm\mu}^{k}+ \beta_{k}({\bm\mu}^{k}-{\bm\mu}^{k-1}).
\end{eqnarray*}
\item[\bf Step 3] If a termination criterion is not met, set $k:=k+1$ and go to Step 1
\end{description}
\end{algorithm}

As we know, {it is important to appropriately choose the two operators $D_{1h}$ and $D_{2h}$ for
efficient numerical computation}. Firstly, if we choose $D_{2h}=0$, which is similar to choosing $D_2=0$
{for the continuous problem in the previous section}, it is unfortunately
{not a  good choice since} there does not exist a closed form solution for the $\bm\mu$-subproblem because the mass matrix $M_h$ is not diagonal. In order to make the {$\bm\mu$-subproblem to have an} analytical solution, we choose
\begin{equation*}
  D_{2h}:=\frac{1}{\alpha}\gamma M_hW_h^{-1}M_h-\frac{1}{\alpha}M_h,\quad {\rm where} \  \gamma =
 \left\{ \begin{aligned}
         &4  \quad if \ n=2, \\
         &5  \quad if \ n=3.
                           \end{aligned} \right.
\end{equation*}
From Proposition \ref{eqn:martix properties1}, it is easy to see that $\mathcal{D}_{2h}\succeq0$.
Let us denote ${\bm\xi}=M_h{\bm\mu}$, then solving the subproblem {for the variable $\bm\mu$} can be translated to solving the following subproblem:
\begin{equation}\label{subproblem-z}
  \begin{aligned}
{\bm\xi}^{k}
     &=\arg\min\frac{1}{2\alpha}\|{\bm\xi}-M_h({\bf p}^k-{\bm\lambda}^k)\|_{M_h^{-1}}^2+\delta^*_{[a,b]}(\bm\xi)+\frac{1}{2\alpha}\|{\bm\xi}-\tilde{\bm\xi}^k\|^2_{\gamma W_h^{-1}-M_h^{-1}}\\
     &=\arg\min\frac{1}{2\alpha}\|{\bm\xi}-(\tilde{\bm\xi}^k+\frac{1}{\gamma}W_h({\bf p}^k-{\bm \lambda}^k-M_h^{-1}\tilde{{\bm\xi}}^k))\|_{\gamma W_h^{-1}}^2+\delta^*_{[a,b]}(\bm\xi).\\
     &={\bf v}^k-\frac{\alpha}{\gamma}W_h{\rm\Pi}_{[a,b]}(\frac{\gamma}{\alpha}W_h^{-1}{\bf v}^k).
  \end{aligned}
\end{equation}
where
\begin{equation*}
  {\bf v}^k= M_h\tilde{\bm\mu}^k+\frac{1}{\gamma}W_h({\bf p}^k-{\bm\lambda}^k-\tilde{\bm\mu}^k).
\end{equation*}
Then we can compute ${\bm\mu}^k$ by ${\bm\mu}^k= M_h^{-1}{\bm\xi}^{k}$.

Next, we {discuss} how to choose the operator $D_{1h}$. {Similar to \eqref{sGS subproblem for D}}, the $({\bm\lambda},{\bf p})$-subproblem can also be rewritten in the following form:
\begin{equation}\label{sGS subproblem for Dh}
  \min \delta_{[-\beta,\beta]}(\bm\lambda)+\frac{1}{2}\Big\langle \left(\begin{array}{c}
                                                             {\bm\lambda} \\
                                                             {\bf p}
                                                           \end{array}\right)
  ,\mathcal{H}_h\left(\begin{array}{c}
                                                             {\bm\lambda}\\
                                                             {\bf p}
                                                           \end{array}\right)\Big\rangle-\Big\langle {\bf r},\left(\begin{array}{c}
                                                             {\bm\lambda} \\
                                                            {\bf p}
                                                           \end{array}\right)\Big\rangle,
\end{equation}
where $\mathcal{H}_h=
\frac{1}{\alpha}\left(
\begin{array}{cc}
M_h& \quad-M_h\\
-M_h & \quad M_h+\alpha K_h M_h^{-1}K_h
\end{array}
\right)$ and ${\bf r}=\left(\begin{array}{c}
                                                             -\frac{1}{\alpha}M_h\tilde{\bm\mu}^k \\
                                                             -M_h{\bf y_r}+K_h{\bf y_d}+\frac{1}{\alpha}M_h\tilde{\bm\mu}^k
                                                           \end{array}\right)$.
Based on the structure of the $({\bm\lambda},{\bf p})$-subproblem, we also use the block sGS decomposition technique to solve it. Thus, we choose
\begin{equation*}
  \mathcal{\widetilde{D}}_{1h}=\mathrm{sGS}(\mathcal{H}_h)=\frac{1}{\alpha}\left(
                                      \begin{array}{cc}
                                        M_h(M_h+\alpha K_h M_h^{-1}K_h)^{-1}M_h &  \quad0\\
                                        0 &  \quad 0\\
                                      \end{array}
                                    \right).
\end{equation*}
And once again, according to \cite[Theorem 2.1]{SunToh3}, we can solve the $({\bm\lambda},{\bf p})$-subproblem by the following steps:
\begin{equation*}
\left\{\begin{aligned}
 \hat{\bf p}^{k}&=\arg\min\frac{1}{2}\|K_h {\bf p}-{M_h}{\bf y_{d}}\|_{M_h^{-1}}^2+ \frac{1}{2\alpha}\|{\bf p}-\tilde{\bm\lambda}^k-\tilde{\bm\mu}^k\|_{M_h}^2+\langle M_h {\bf y_r, p}\rangle,\\
 {\bm\lambda}^{k}
        &=\arg\min\frac{1}{2\alpha}\|{\bm\lambda}-(\hat{\bf p}^{k}-\tilde{\bm\mu}^k)\|_{M_h}^2+\delta_{[-\beta,\beta]}(\bm\lambda),\\
{\bf p}^{k}&=\arg\min\frac{1}{2}\|K_h {\bf p}-{M_h}{\bf y_{d}}\|_{M_h^{-1}}^2+ \frac{1}{2\alpha}\|{\bf p}-{\bm\lambda}^k-\tilde{\bm\mu}^k\|_{M_h}^2+\langle M_h {\bf y_r, p}\rangle.
\end{aligned}\right.
\end{equation*}

However, it is easy to see that the $\bm\lambda$-subproblem is {not a simple projection problem} with respect to the variable $\bm\lambda$ since the mass matrix $M_h$ is not diagonal, thus there is no closed form solution for $\bm\lambda$. To overcome this difficulty, we can add a proximal term $\frac{1}{2\alpha}\|\bm\lambda- \tilde{\bm\lambda}^{k}\|_{W_h-M_h}^2$ to the $\bm\lambda$-subproblem.
{Then} for the $\bm\lambda$-subproblem, we have
\begin{equation*}
  {\bm\lambda}^{k}={\rm\Pi}_{[-\beta,\beta]}(\tilde{\bm\lambda}^{k}+W_h^{-1}M_h(\hat{\bf p}^k-\tilde{\bm\mu}^k-\tilde{\bm\lambda}^{k})).
\end{equation*}
Thus, we can choose $D_{1h}$ as follows
\begin{equation*}
    \mathcal{D}_{1h}={\mathrm{sGS}\left(\mathcal{H}_h +
    \frac{1}{\alpha}\left[
                                      \begin{array}{cc}
                                        W_h-M_h &\quad0  \\
                                        0 & \quad0 \\
                                      \end{array}
                                    \right]\right) +\left( \frac{1}{\alpha}\left[
                                      \begin{array}{cc}
                                        W_h-M_h &\quad0  \\
                                        0 & \quad0 \\
                                      \end{array}
                                    \right]\right)}
\end{equation*}
Then, according to the  {above choices} of $\mathcal{D}_{1h}$ and $\mathcal{D}_{2h}$,  the detailed framework of our inexact sGS based majorized ABCD method for (\ref{eqn:discretized matrix-vector dual problem}) {is given} as follows.
\begin{algorithm}[H]
  \caption{\textbf{(sGS-mABCD algorithm for (\ref{eqn:discretized matrix-vector dual problem}))}}\label{algo1:Full inexact ABCD algorithm for (Dh)}
  \textbf{Input}:{$(\tilde{\bm\lambda}^1, \tilde{\bf p}^1,\tilde{\bm\mu}^1)=({\bm\lambda}^0, {\bf p}^0,\bm\mu^0)\in {\rm dom} (\delta^*_{[a,b]})\times [-\beta,\beta]\times \mathbb{R}^{N_h}$. Set $k= 1, t_1= 1.$}

 \textbf{Output}:{$ ({\bm\lambda}^k, {\bf p}^k,{\bm\mu}^k)$}
\begin{description}
  \item[Step 1] Compute 
\begin{eqnarray*}
 \hat{\bf p}^{k}&=&\arg\min\frac{1}{2}\|K_h {\bf p}-{M_h}{\bf y_{d}}\|_{M_h^{-1}}^2+ \frac{1}{2\alpha}\|{\bf p}-\tilde{\bm\lambda}^k-\tilde{\bm\mu}^k\|_{M_h}^2+\langle M_h {\bf y_r, p}\rangle,\\
 {\bm\lambda}^{k}
        &=&\arg\min\delta_{[-\beta,\beta]}(\bm\lambda)+\frac{1}{2\alpha}\|\bm\lambda-(\hat{\bf p}^{k}-\tilde{\bm\mu}^k)\|_{M_h}^2+\frac{1}{2\alpha}\|\bm\lambda- \tilde{\bm\lambda}^{k}\|_{W_h-M_h}^2,\\
{\bf p}^{k}&=&\arg\min\frac{1}{2}\|K_h {\bf p}-{M_h}{\bf y_{d}}\|_{M_h^{-1}}^2+ \frac{1}{2\alpha}\|{\bf p}-{\bm\lambda}^k-\tilde{\bm\mu}^k\|_{M_h}^2+\langle M_h {\bf y_r, p}\rangle,\\
  {\bm \mu}^{k}&=&\arg\min\delta^*_{[a,b]}(M_h\bm\mu)+\frac{1}{2\alpha}\|\bm\mu-({\bf p}^k-\bm\lambda^k)\|_{M_h}^2+\frac{1}{2\alpha}\|\bm\mu-\tilde{\bm\mu}^k\|^2_{\gamma M_hW_h^{-1}M_h-M_h}.
\end{eqnarray*}

  \item[Step 2] Set $t_{k+1}=\frac{1+\sqrt{1+4t_k^2}}{2}$ and $\beta_k=\frac{t_k-1}{t_{k+1}}$, compute
\begin{eqnarray*}
\tilde{\bm\lambda}^{k+1}= {\bm\lambda}^{k}+ \beta_{k}({\bm\lambda}^{k}-{\bm\lambda}^{k-1}),\quad
\tilde{\bf p}^{k+1}={\bf p}^{k}+ \beta_{k}({\bf p}^{k}-{\bf p}^{k-1}), \quad\tilde{\bm\mu}^{k+1}={\bm\mu}^{k}+ \beta_{k}({\bm\mu}^{k}-{\bm\mu}^{k-1}).
\end{eqnarray*}
\item[\bf Step 3] If a termination criterion is not met, set $k:=k+1$ and go to Step 1
\end{description}
\end{algorithm}
Similarly, owing to Theorem \ref{imABCD convergence}, we can show Algorithm \ref{algo1:Full inexact ABCD algorithm for (Dh)} also has the following $O(1/k^2)$ iteration complexity.
\begin{theorem}\label{sGS-imABCD convergence}
Suppose that the solution set $\Omega$ of the problem {\rm(\ref{eqn:discretized matrix-vector dual problem})} is non-empty. Let
${\bf z}^*=(\bm\lambda^*,{\bf p}^*,{\bm\mu}^*)\in \Omega$. Let $\{{\bf z}^k\}:=\{({\bm\lambda}^k,{\bf p}^k,{\bm\mu}^k)\}$ be the sequence generated by the Algorithm {\rm\ref{algo1:Full inexact ABCD algorithm for (Dh)}}. Then we have
\begin{equation}\label{iteration complexity of discretized dual problem}
\Phi_h({\bf z}^k)- \Phi_h({\bf z}^*)\leq \frac{4\tau_h}{(k+1)^2}, \; \forall k\geq 1,
\end{equation}
where $\Phi_h(\cdot)$ is the objective function of the dual problem {\rm(\ref{eqn:discretized matrix-vector dual problem})} and
\begin{eqnarray}
    &&\tau_h=\frac{1}{2}\|{\bf z}^0- {\bf z}^*\|_{\mathcal{S}_h}^2 \label{tauh}\\
    &&\mathcal{S}_h:=\frac{1}{\alpha}\left(
                                      \begin{array}{ccc}
                                        M_h(M_h+\alpha K_h M_h^{-1}K_h)^{-1}M_h+W_h-M_h &  \ 0&\  0\\
                                        0 &  \ 0&  \ 0\\
                                        0&  \ 0&  \ \gamma M_hW_h^{-1}M_h\\
                                      \end{array}
                                    \right)\label{Sh}
\end{eqnarray}
Moreover, the sequence $\{({\bm\lambda}^k,{\bf p}^k,{\bm\mu}^k)\}$ generated by the Algorithm {\rm\ref{algo1:Full inexact ABCD algorithm for (Dh)}} is bounded.
\end{theorem}

\section{Robustness with respect to $h$}\label{sec:4}
In this section,
{we deal with the issue on} how measures of the convergence behavior of the iteration sequence vary with the level of approximation. Such questions come under the category of mesh-independence results. In this section, we will establish the mesh independence of majorized accelerate block coordinate descent (mABCD) method for optimal control problems.

In what follows we will give one type of mesh-independence result for mABCD method.
{
It says that the iterate $k$ after which the difference $\Phi_h({\bf z}^k)-\inf\Phi_h({\bf z})$ has been identified up to less than $\epsilon$ is independent of the mesh size $h$.
}
In order to show these results, let us first present some bounds on the Rayleigh quotients of $K_h$ and $M_h$, one can see \cite[Proposition 1.29 and Theorem 1.32]{spectralproperty} for more details.

\begin{lemma}\label{spectral property}
For $\mathcal{P}1$ approximation on a regular and quasi-uniform subdivision of $\mathbb{R}^n$ which satisfies Assumption {\rm\ref{regular and quasi-uniform triangulations}}, and for any ${\bf x}\in \mathbb{R}^{N_h}$, the mass matrix $M_h$ approximates the scaled identity matrix in the sense that
\begin{equation*}
c_1 h^2\leq \frac{{\bf x}^{T}M_h{\bf x}}{{\bf x}^{T}{\bf x}}\leq c_2 h^2 \quad if\  n=2, \ {\rm and}\ c_1 h^3\leq \frac{{\bf x}^{T}M_h{\bf x}}{{\bf x}^{T}{\bf x}}\leq c_2 h^3 \quad if\  n=3,
\end{equation*}
the stiffness matrix $K_h$ satisfies
\begin{equation*}
d_1h^2\leq \frac{{\bf x}^{T}K_h{\bf x}}{{\bf x}^{T}{\bf x}}\leq d_2 \quad if\  n=2, \ {\rm and }\ d_1h^3\leq \frac{{\bf x}^{T}K_h{\bf x}}{{\bf x}^{T}{\bf x}}\leq d_2 h \quad if\  n=3,
\end{equation*}
where the constants $c_1$, $c_2$, $d_1$ and $d_2$ are independent of the mesh size $h$.
\end{lemma}

Based on Lemma \ref{spectral property}, we can easily obtain the following lemma.
\begin{lemma}\label{spectral property of Gh}
Let $G_h=M_h+\alpha K_h M_h^{-1}K_h$. For any ${\bf x}\in \mathbb{R}^{N_h}$, there exist four constants $u_1$, $u_2$, $l_1$, $l_2$ and $h_0>0$, such that for any $0<h<h_0$, the matrix $G_h$ satisfies the following inequalities
\begin{equation}\label{spectral property of Gh12}
  \begin{aligned}
&l_1 h^2\leq \frac{{\bf x}^{T}G_h{\bf x}}{{\bf x}^{T}{\bf x}}\leq u_1 \frac{1}{h^2} \quad if\  n=2, \quad
l_2 h^3\leq \frac{{\bf x}^{T}G_h{\bf x}}{{\bf x}^{T}{\bf x}}\leq u_2 \frac{1}{h} \quad\ if\  n=3.
  \end{aligned}
\end{equation}
\end{lemma}

Thus based on Lemma \ref{spectral property} and Lemma \ref{spectral property of Gh}, it is easy to prove that there exists $h_0>0$, such that for any $0<h<h_0$, the matrix $M_hG_h^{-1}M_h$ satisfies the following properties
\begin{equation*}\label{property of MinvGM}
\begin{aligned}
 &\lambda_{\max}(M_hG_h^{-1}M_h)=O(h^2) \quad {\rm for}\ n=2,\quad \lambda_{\max}(M_hG_h^{-1}M_h)=O(h^3) \quad {\rm for}\ n=3,
\end{aligned}
\end{equation*}
where $ \lambda_{\max}(\cdot)$ {denotes} the largest eigenvalue of a given matrix. Furthermore, we have
\begin{equation}\label{property of Sh}
\begin{aligned}
  \lambda_{\max}(\mathcal{S}_h)&=\frac{1}{\alpha}\max\{ \lambda_{\max}(M_hG_h^{-1}M_h+W_h-M_h),\lambda_{\max}(\gamma M_hW_h^{-1}M_h)\}\\
  &=\left\{\begin{aligned}
      & O(h^2) \quad {\rm for}\ n=2,\\
      & O(h^3) \quad {\rm for}\ n=3.
    \end{aligned}\right.
\end{aligned}.
\end{equation}
where $\mathcal{S}_h$ defined in (\ref{Sh}). {In other words, we can say that the largest eigenvalue of the matrix $\mathcal{S}_h$ can be uniformly bounded by a constant},
which implies the ``discretized" convergence factor $\tau_h$
could be uniformly bounded by a constant. Hence, this conclusion prompts us to consider analysing the mesh independence of {the}
mABCD method. We present our first mesh independence result for our mABCD method, in which we prove that the ``discretized" convergence factor $\tau_h$ defined in Theorem \ref{sGS-imABCD convergence} approach the ``continuous" convergence factor $\tau$ defined in {Theorem \ref{imABCD convergence for dual problem}}
 in the limits $h\rightarrow 0$ and the distance can be bounded in terms of the mesh size.

\begin{theorem}\label{mesh independent result1}
Let Algorithm {\rm\ref{algo1:ABCD algorithm for (D)}} for {the} continuous problem {\rm(\ref{eqn:dual problem})} start from $z^0=({\lambda}^0, {p}^0,\mu^0)$ and Algorithm {\rm\ref{algo1:Full inexact ABCD algorithm for (Dh)}} for {the} discretized problem {\rm(\ref{eqn:discretized matrix-vector dual problem})} start from $z_h^0=({\lambda}_h^0, {p}_h^0,\mu_h^0)$, respectively. And we take $z(x)^* \in(\partial\Phi)^{-1}(0)$ and $z_h^*(x)=\sum\limits_{i=1}^{N_h}z_i^*\phi(x)$ where the coefficients $(z_1^*, z_2^*,...,z_{N_h}^*)\in(\partial\Phi_h)^{-1}(0)$. Assume that $z_h^0=I_hz^0$ where $I_h$ is the nodal interpolation operator, and $\|z^*-z_h^*\|_{L^2(\Omega)}=O(h)$. Then there exist $h^*\in (0, \hat h]$ and a constant $C$, such that
\begin{equation}\label{difference of tauh and tau}
  \tau_h\leq\tau+Ch
\end{equation}
for all $h \in(0, h^*]$.
\begin{proof}
From the definition of $\tau$ in Theorem \ref{imABCD convergence}, we have
\begin{equation}\label{defintion of tau}
\begin{aligned}
  \tau&=\frac{1}{2\alpha}\|\mu^*-\mu^0\|^2_{L^2(\Omega)}+\frac{1}{2\alpha}\langle\lambda^*-\lambda^0, (\alpha A^*A+\mathcal{I})^{-1}(\lambda^*-\lambda^0)\rangle_{L^2(\Omega)}\\
      &=\frac{1}{2\alpha}\|\mu^*-\mu^0\|^2_{L^2(\Omega)}+\frac{1}{2\alpha}\int_{\Omega}(\lambda^*-\lambda^0)q^1{~\rm dx},
\end{aligned}
\end{equation}
where $q^1$ is the weak solution of the following problem:
\begin{eqnarray}\label{PDEs1}
\nonumber&&{\rm Find}~(q^1,q^2)\in (H^1_0(\Omega))^2, {\rm such~ that}\\
 &&\left\{\begin{aligned}
    &a(q^1, v) =\langle q^2, v\rangle_{L^2(\Omega)},\\
    &\alpha a(q^2, v)+\langle q^1, v\rangle_{L^2(\Omega)}=\langle \lambda^*-\lambda^0, v\rangle_{L^2(\Omega)}, \forall v\in H^1_0(\Omega)
  \end{aligned}\right.
\end{eqnarray}
where the bilinear form $a(\cdot,\cdot)$ is defined in (\ref{eqn:bilinear form}). Similarly, according to the definition of $\tau_h$ and Proposition \ref{eqn:martix properties}, we obtain
\begin{equation}\label{defintion of tauh}
\begin{aligned}
  \tau_h&=\frac{1}{2\alpha}\|\mu_h^*-\mu_h^0\|^2_{L^2(\Omega)}+\frac{1}{2\alpha}\int_{\Omega}(\lambda_h^*-\lambda_h^0)q^1_h{~\rm dx}+\frac{1}{2\alpha}\|I_h(\lambda_h^*-\lambda_h^0)^2-(\lambda_h^*-\lambda_h^0)^2\|_{L^1(\Omega)},
\end{aligned}
\end{equation}
where $q_h^1$ is the solution of the following discretized problem which is discretized by piecewise linear finite elements:
\begin{eqnarray}\label{PDEs2}
\nonumber&&{\rm Find}~(q_h^1,q_h^2), {\rm such~ that}\\
 &&\left\{\begin{aligned}
    &a(q_h^1, v_h) =\langle q_h^2, v_h\rangle_{L^2(\Omega)},\\
    &\alpha a(q_h^2, v_h)+\langle q_h^1, v_h\rangle_{L^2(\Omega)}=\langle \lambda_h^*-\lambda_h^0, v_h\rangle_{L^2(\Omega)}, \forall v_h\in  Y_h
  \end{aligned}\right.
\end{eqnarray}
where $Y_h$ is defined in {\rm(\ref{state discretized space})}. In order to estimate the value of $\tau_h$, we define
$\tilde{q}_h^1$ as the solution of the following discretized problem
\begin{eqnarray}\label{PDEs3}
\nonumber&&{\rm Find}~(\tilde{q}_h^1,\tilde{q}_h^2), {\rm such~ that}\\
 &&\left\{\begin{aligned}
    &a(\tilde{q}_h^1, v_h) =\langle \tilde{q}_h^2, v_h\rangle_{L^2(\Omega)},\\
    &\alpha a(\tilde{q}_h^2, v_h)+\langle \tilde{q}_h^1, v_h\rangle_{L^2(\Omega)}=\langle \lambda^*-\lambda^0, v_h\rangle_{L^2(\Omega)}, \forall v_h\in  Y_h.
  \end{aligned}\right.
\end{eqnarray}
Obviously, there exists $h^*\in (0, \hat h]$ and four constants $C_1$ and $C_2$, $C_3$ and $C_4$ which independent of $h$, such that for all $h \in(0, h^*]$, the following inequalities hold:
\begin{equation}\label{triangle inequality}
\begin{aligned}
  \|q^1-q_h^1\|_{L^2(\Omega)}&\leq\|q^1-\tilde{q}_h^1\|_{L^2(\Omega)}+\|q_h^1-\tilde{q}_h^1\|_{L^2(\Omega)}\\
  &\leq C_1h^2\|\lambda^*-\lambda^0\|_{L^2(\Omega)}+C_2h^2\|\lambda_h^*-\lambda_h^0\|_{L^2(\Omega)}\\
  &\leq C_1h^2\|\lambda^*-\lambda^0\|_{L^2(\Omega)}+C_2h^2(\|\lambda_h^*-\lambda^*\|_{L^2(\Omega)}+\|\lambda_h^0-\lambda^0\|_{L^2(\Omega)}+\|\lambda^*-\lambda^0\|_{L^2(\Omega)})\\
  &\leq C_3h^2\|\lambda^*-\lambda^0\|_{L^2(\Omega)}+C_4h^3(\|\lambda^*\|_{L^2(\Omega)}+\|\lambda^0\|_{L^2(\Omega)}).
\end{aligned}
\end{equation}
Thus, we now can estimate $\tau$ and get
\begin{equation}\label{estimate tau}\small
\begin{aligned}
\tau_h=&\frac{1}{2\alpha}\|\mu_h^*-\mu^*+\mu^0-\mu_h^0+\mu^*-\mu^0\|^2_{L^2(\Omega)}+\frac{1}{2\alpha}\int_{\Omega}(\lambda^*-\lambda^0)q^1{~\rm dx}+\frac{1}{2\alpha}\int_{\Omega}(\lambda_h^*-\lambda^*)(q_h^1-q^1){~\rm dx}\\&+\frac{1}{2\alpha}\int_{\Omega}(\lambda^0-\lambda_h^0)(q_h^1-q^1){~\rm dx}+\frac{1}{2\alpha}\int_{\Omega}(\lambda^*-\lambda^0)(q_h^1-q^1){~\rm dx}+\frac{1}{2\alpha}\int_{\Omega}(\lambda^0-\lambda_h^0)q^1{~\rm dx}\\
&+\frac{1}{2\alpha}\int_{\Omega}(\lambda_h^*-\lambda^*)q^1{~\rm dx}+\frac{1}{2\alpha}\|I_h(\lambda_h^*-\lambda_h^0)^2-(\lambda_h^*-\lambda_h^0)^2\|_{L^1(\Omega)}\\
\leq&\frac{1}{2\alpha}\|\mu^*-\mu^0\|^2_{L^2(\Omega)}+\frac{1}{2\alpha}\int_{\Omega}(\lambda^*-\lambda^0)q^1{~\rm dx}+\frac{1}{2\alpha}\|\mu_h^*-\mu^*\|^2_{L^2(\Omega)}+\frac{1}{2\alpha}\|\mu^0-\mu_h^0\|^2_{L^2(\Omega)}\\
&+\frac{1}{2\alpha}\int_{\Omega}(\lambda_h^*-\lambda^*)(q_h^1-q^1){~\rm dx}
+\frac{1}{2\alpha}\int_{\Omega}(\lambda^0-\lambda_h^0)(q_h^1-q^1){~\rm dx}
+\frac{1}{2\alpha}\int_{\Omega}(\lambda^*-\lambda^0)(q_h^1-q^1){~\rm dx}\\
&+\frac{1}{2\alpha}\int_{\Omega}(\lambda^0-\lambda_h^0)q^1{~\rm dx}
+\frac{1}{2\alpha}\int_{\Omega}(\lambda_h^*-\lambda^*)q^1{~\rm dx}
+\frac{1}{2\alpha}\|I_h(\lambda_h^*-\lambda_h^0)^2-(\lambda_h^*-\lambda_h^0)^2\|_{L^1(\Omega)}\\
\leq&\tau+C_5h(\|\lambda^0\|_{L^2(\Omega)}+\|\lambda^*\|_{L^2(\Omega)})+C_6h^2(\|\mu^*\|_{L^2(\Omega)}+\|\mu^0\|_{L^2(\Omega)}+\|\lambda^0\|_{L^2(\Omega)}+\|\lambda^*\|_{L^2(\Omega)})+O(h^3)\\
\leq&\tau+Ch
\end{aligned}
\end{equation}
\end{proof}
\end{theorem}

%
\section{Concluding remarks}\label{sec:5}
In this paper, instead of solving the optimal control problem with $L^1$ control cost, we directly solve {its dual, which is a
multi-block unconstrained convex composite minimization problem}. By taking advantage of the structure of the dual problem, and combining the majorized ABCD (mABCD) method and the recent advances in the inexact symmetric Gauss-Seidel (sGS) technique, we introduce the sGS-mABCD method to solve the dual problem. More importantly, one type of mesh independence {result for the}
mABCD method is proved, which asserts that asymptotically the infinite dimensional {mABCD} method and
{the} finite dimensional discretization {version} have the same convergence property
{in the sense that the worst case iteration complexity of the mABCD method remains} nearly constant as the discretization is refined.

\end{document}